\documentclass[12pt]{amsart}

\usepackage[headings]{fullpage}
\usepackage[utf8]{inputenc}
\usepackage{hyperref}
\usepackage{url}
\usepackage{mathtools}
\usepackage{stmaryrd}
\usepackage{amsmath,amssymb,amsthm,amscd}

\title[Bounded functions on the character variety]{Bounded functions on the character variety for quadratic extensions of $\mathbb{Q}_p$}

\author{Laurent Berger}
\address{UMPA, ENS de Lyon, 
UMR 5669 du CNRS, 69007 Lyon}
\email{laurent.berger@ens-lyon.fr}

\author{Johannes Sprang}
\address{Universit\"at Duisburg-Essen, 
Fakult\"at f\"ur Mathematik, 
Thea-Leymann-Str. 9, D-45127 Essen}
\email{johannes.sprang@uni-due.de}

\author{Rustam Steingart}
\address{Ruprecht-Karls-Universit\"at Heidelberg,
Mathematisches Institut,
Im Neuenheimer Feld 205, D-69120 Heidelberg}
\email{rsteingart@mathi.uni-heidelberg.de}

\date{\today}

\newtheorem{theorem}{Theorem}[section]
\newtheorem*{theorem*}{Theorem}
\newtheorem{lemma}[theorem]{Lemma}

\newtheorem{corollary}[theorem]{Corollary}

\theoremstyle{remark}

\newtheorem{question}[theorem]{\bf Question}
\newtheorem{conjecture}[theorem]{\bf Conjecture}

\newcommand{\Cp}{\mathbb{C}_p}
\newcommand{\Qp}{\mathbb{Q}_p}
\newcommand{\Zp}{\mathbb{Z}_p}
\newcommand{\ol}[1]{\overline{#1}}

\newcommand{\spc}{\mathrm{sp}}
\newcommand{\quot}{\mathrm{q}}

\newcommand{\ul}[1]{\underline{#1}}

\DeclareMathOperator{\Sp}{Sp}
\DeclareMathOperator{\gr}{gr}
\DeclareMathOperator{\im}{im}
\DeclareMathOperator{\End}{End}
\newcommand{\frX}{\mathfrak{X}}
\newcommand{\frY}{\mathfrak{Y}}

\newcommand{\dcroc}[1]{[\![ #1 ]\!]}
\newcommand{\NN}{\mathbb{N}}

\newcommand{\Gal}{\operatorname{Gal}}
\newcommand{\Hom}{\operatorname{Hom}}
\newcommand{\ZZ}{\mathbb{Z}}

\newcommand{\QQ}{\mathbb{Q}}

\newcommand{\RR}{\mathbb{R}}

\newcommand{\calA}{\mathcal{A}}
\newcommand{\calM}{\mathcal{M}}
\newcommand{\calQ}{\mathcal{Q}}

\newcommand{\Iwalg}{o_L\llbracket o_L\rrbracket}

\begin{document}

\begin{abstract}
For any finite extension $L$ of $\Qp$, Schneider--Teitelbaum defined and studied the character variety in $p$-adic Fourier theory. It is a rigid analytic curve over $L$ that parameterizes the set of locally $L$-analytic characters of $o_L$. It is conjectured that all bounded functions on the character variety of Schneider--Teitelbaum come from $p$-adic measures. We prove this conjecture for all quadratic extensions $L/\Qp$.
\end{abstract}

\maketitle 

\section{Introduction}
Let $L$ be a finite extension of $\Qp$ and let $\Cp$ denote the completion of an algebraic closure of $\Qp$. In their work on $p$-adic Fourier theory \cite{SchneiderTeitelbaum2001}, Schneider and Teitelbaum defined and studied the character variety $\frX$. This character variety is a rigid analytic curve over $L$ that parameterizes the set of locally $L$-analytic characters $\lambda : (o_L,+) \to (\Cp^\times,\times)$. 

The character variety plays an important role in the theory of $p$-adic representations and in the construction of $p$-adic $L$-functions. Schneider and Teitelbaum show that the space of functions $\mathcal{O}(\frX)$ on the character variety identifies with the space of locally $L$-analytic distributions $D(o_L,L)$. However, we have so far a fairly poor understanding of the bounded-by-$1$ functions $\mathcal{O}(\frX)^\circ$ on the character variety. Since every $p$-adic measure gives rise to a unique bounded locally $L$-analytic distribution, we have an injective map of $o_L$-algebras
\begin{equation}
o_L\llbracket o_L \rrbracket \hookrightarrow \mathcal{O}(\frX)^\circ.
\end{equation}
A natural question first raised by Teitelbaum and studied extensively in \cite{ArdakovBerger2024} is the following:

\begin{question} \label{question:main}
Do we have $o_L\llbracket o_L\rrbracket=\mathcal{O}(\frX)^\circ$?
\end{question}

Some numerical evidence for a positive answer to this question has been given in \cite{BergerSprang2026} for quadratic unramified extensions. In this paper, we show that Question \ref{question:main} has a positive answer for any quadratic extension $L$ of $\Qp$. More precisely:

\begin{theorem}\label{thm:main}
If $L/\Qp$ is a quadratic extension, the canonical injective map $o_L\llbracket o_L\rrbracket\hookrightarrow \mathcal{O}(\frX)^\circ$ is an isomorphism.
\end{theorem}

Before we discuss some applications of Theorem \ref{thm:main}, let us briefly summarize the key ideas of its proof. For simplicity, assume that $L$ is a quadratic unramified extension. The proof in the ramified case is only a minor variation of the argument in the unramified case. In the quadratic unramified case, the locally $L$-analytic character variety can be cut out from the locally $\Qp$-analytic character variety  $\mathbb{D}^\circ \otimes_{\Zp}\Hom(o_L,\Zp)\cong (\mathbb{D}^\circ)^2$ by a single equation $H$. On the level of affinoid polydiscs of radius $\rho_n\in (0,1)$, this gives a surjection
\[
	B_n=L\langle \rho_n^{-1}x_1,\rho_n^{-1}x_2 \rangle\twoheadrightarrow A_n:=B_n/(H).
\] 
The key observation of this paper is that for suitably chosen radii $\rho_n$, we get a surjection $B_n^\circ \twoheadrightarrow A_n^\circ$ for the rings of bounded functions. In other words, we can lift any bounded function in $A_n$ to a bounded function in $B_n$. We prove this by comparing the quotient norm on $A_n$ induced by $B_n\twoheadrightarrow A_n$ to the spectral norm on $A_n$. Using this, Theorem \ref{thm:main} is essentially a consequence of the isomorphism $\mathcal{O}(\frX)^\circ \cong \varprojlim_{(0,1)\ni \rho_n\to 1} A_n^\circ$ and the fact that we can describe the bounded functions $B_n^\circ$ on a polydisc of radius $\rho_n$ explicitly. In section \ref{sec:concluding-remarks}, we explain the limitations of our method. In particular, we show that the quotient norm on $A_n$ is in general not the spectral norm of $A_n$.

In \cite{ArdakovBerger2024}, Ardakov and the first-named author have established several connections between Question \ref{question:main} and Katz's theory of Galois-continuous measures as well as with questions about integer-valued polynomials. For example, they show that a positive answer to Question \ref{question:main} implies the following conjecture about integer-valued polynomials:

\begin{conjecture}\label{conjecture:intpol}
For any finite extension $L/\Qp$ we have $\mathbf{Int}=\mathbf{Pol}$.
\end{conjecture}

Here, $\mathbf{Int} = \{ P(X) \in L[X]$ such that $P(o_L) \subset o_L \}$ is the space of all integer-valued polynomials on $o_L$, and $\mathbf{Pol}$ denotes the $o_L$-module spanned by the polynomials $c_{n,i}(X)$ for all $i \geq n \geq 0$, where $[a](T)^n = \sum_{i \geq n} c_{n,i}(a) T^i$, where $[a](T)$ denotes the multiplication-by-$a$-map on a Lubin--Tate formal $o_L$-module attached to a uniformizer of $o_L$. Hence, our main theorem also proves Conjecture \ref{conjecture:intpol} for all quadratic extensions of $\Qp$. Note that a stronger implication is proved in remark 3.5 of \cite{BergerSprang2026}.

Another implication of Theorem \ref{thm:main} concerns associated orders in a Lubin--Tate extension. Let $L = \QQ_{p^2}$ and let $L_n$ be the extension of $L$ generated by the $p^n$-torsion points of the Lubin--Tate group attached to $L$ and $p$. Let $G_n = \Gal(L_n/L)$ and let $\calA_n = \{ f \in L[G_n]$ such that $f(o_n) \subset o_n \}$. There are restriction maps $\pi_n : \calA_{n+m} \to \calA_n$ for all $m,n \geq 0$. In \S \ref{assord}, we prove the following result.

\begin{theorem}
\label{ordcrit:intro}
For all $n \geq 0$, there exists $m=m(n)$ such that $\pi_n(\calA_{n+m}) = o_L[G_n]$.
\end{theorem}

At about the same time that we proved Theorem \ref{thm:main}, we learned that Uro\v{s} Colovi\'c had answered Question \ref{question:main} for all finite extensions $L/\Qp$ using a completely different method.

\section{Notation and Preliminaries}\label{sec:notation}
Let $L/\QQ_p$ be a finite extension of degree $d=ef$, where $f$ denotes the inertia degree and $e$ the ramification index of $L/\QQ_p$. The field $L$ has valuation ring $o_L$, uniformiser $\pi$, and residue field $k_L$. We write $q:=|k_L|=p^f$ for the cardinality of the residue field. We denote by $\Iwalg$ the completed group ring of $o_L$. Fixing a basis $(v_1,\dots,v_d)$ of $o_L$ over $\ZZ_p$, we obtain an  isomorphism $o_L\llbracket x_1,\dots,x_d\rrbracket \cong \Iwalg$, by sending $x_i$ to $\delta_{v_i}-1$, where $\delta_{x}$ denotes the Dirac measure which evaluates a continuous function at $x$. We denote by $\mathbb{D}^\circ$ the rigid analytic open unit disc over $L$. 

From the above isomorphism we can see that the rigid analytic generic fibre $\frY$ of $\operatorname{Spf}(\Iwalg)$ is isomorphic to the open unit polydisc $(\mathbb{D}^\circ)^d$. For a complete field extension $C$ of $L$, the $C$-valued points of $\frY$ parametrize locally $\QQ_p$-analytic characters $o_L \to C^\times$.  In \cite{SchneiderTeitelbaum2001}, the authors construct a closed rigid analytic subvariety $\frX \subseteq \frY$ which parametrizes the locally $L$-analytic characters. Since our proof makes use of the explicit description of $\mathcal{O}(\frX)$, we recall the construction of $\frX$ in the following. The ring $\mathcal{O}(\frY)$ is isomorphic to the algebra of locally $\QQ_p$-analytic distributions $D^{\QQ_p}(o_L,L)$. Recall that the Lie-algebra $\mathfrak{g} \cong L$ of $o_L$ can be embedded into $D^{\QQ_p}(o_L,L)$ and a locally $\QQ_p$-analytic character $o_L \to V$ valued in a Banach space over $L$ is locally $L$-analytic if and only if its derivative at $0$ is $L$-linear. Let us write $\iota \colon \mathfrak{g} \to D^{\QQ_p}(o_L,L)$.
By the construction of Schneider--Teitelbaum $\frX$ is the rigid analytic variety corresponding to the ideal $\mathcal{I}$ generated by elements of the form $(t\iota(\mathfrak{x})-\iota(t\mathfrak{x}))$, where $t \in L$ and $\mathfrak{x} \in \mathfrak{g}$.
Because $\mathfrak{g}$ is one-dimensional over $L$ and $(v_1,\dots,v_d)$ is a basis of $L$ over $\Qp$ one has for any non-zero element $\mathfrak{x} \in \mathfrak{g}$
\[
\mathcal{I} = (\iota(v_1\mathfrak{x})-v_1\iota(\mathfrak{x}),\dots,\iota(v_d\mathfrak{x})-v_d\iota(\mathfrak{x})).
\] 
Note that by choosing $v_1:=1$ we can always eliminate the first equation.
By \cite[Lemma 4.3]{SchneiderTeitelbaum2002} the isomorphism
\[ 
D^{\QQ_p}(o_L,L) \cong \mathcal{O}(\frY) \cong \mathcal{O}((\mathbb{D}^\circ)^d)
\]
is compatible with the logarithm map $\log_{o_L} : o_L \to \mathfrak{g}$ in the sense that $\log_{o_L}(v_i)$ is sent to $\log(1+x_i)$.
We denote by $I$ the ideal of $B:=\mathcal{O}((\mathbb{D}^\circ)^d)$ generated by the image of $\mathcal{I}$. By choosing $v_1:=1$ so that $v_j = v_jv_1$ we see that this is precisely the ideal generated by the elements 
$$H_{i}(x_1,\dots,x_d):= \log(1+x_i)-v_i\log(1+x_1).$$

For any increasing sequence $\rho_n$ of elements of $p^{\QQ} \cap (0,1)$ converging to $1$ we can write 
\[
	(\mathbb{D}^\circ )^d=\varinjlim_n \Sp(B_n)
\]
with the Tate-algebras $B_n := L\langle \rho_n^{-1}x_1,\dots,\rho_n^{-1}x_d\rangle$. 
By Theorem B for Stein spaces one has $\mathcal{O}(\frX) \cong \mathcal{O}(\frY)/I \cong \varprojlim B_n/I$.
From this one obtains an explicit exhaustion of $\frX$ by affinoid subdomains $X_n:=\Sp(A_n)$ with $A_n:=B_n/I$.
We can equip $A_n$ with the quotient norm $\lVert \cdot \rVert_{n,\quot}$ induced by the Gau{\ss} norm $\lVert \cdot \rVert_n$ on $B_n$.
Alternatively $A_n$ can be equipped with its spectral semi norm $\lVert \cdot \rVert_{n,\spc}$.
We denote by $$A_n^\circ:=\{f \in A_n \mid \lVert f\rVert_{n,\spc} \leq 1\}$$ the subring of power-bounded elements in $A_n$.

\section{The unramified case}
We specialize to the situation where $L/\Qp$ is a quadratic unramified extension.

Choose a $\Zp$-basis $(v_1=1,v_2)$ of $o_L$. We write $x_1,x_2$ for the coordinates on $\mathbb{D}^\circ\times \mathbb{D}^\circ$, so that the character variety of Schneider--Teitelbaum $\frX$ is locally cut out by the equation
\[
	H(x_1,x_2):=\log(1+x_2)-v_2\log(1+x_1).
\]
For us, it will be convenient to choose the sequence
\[
	\rho_n:=p^{-\frac{1}{p^n-p^{n-1}}}.
\]
of critical radii of the logarithm series, i.e. the increasing sequence $(\rho_n)_n$ consists exactly of those numbers $\rho \in (0,1)$ where the function $\ZZ_{\geq 1} \to \RR, \quad m\mapsto \rho^m/|m|_p$ does not have a unique maximum. Indeed, for $\rho=\rho_n$ the maximum is attained for the two values $m=p^{n-1}$ and $m=p^n$.
As before we have
\[
	\mathbb{D}^\circ\times \mathbb{D}^\circ=\varinjlim_n \Sp(B_n)
\]
with the Tate algebra $B_n=L\langle \rho_n^{-1} x_1, \rho_n^{-1}x_2  \rangle$. The character variety of Schneider and Teitelbaum is given by $\frX=\varinjlim_n \frX_n$ where $\frX_n=\Sp(A_n)$ with
\[
	A_n=B_n/(H).
\]

The goal of this section is to prove that the quotient norm on $A_n$ induced by $B_n \twoheadrightarrow A_n$ is exactly the spectral norm on $A_n$. As a consequence, any bounded function on $A_n$ can be lifted to a bounded function on $B_n$ of the same norm. The computations are very much inspired by \S 3 of Tobias Schmidt's PhD thesis, see \cite{Schmidt2006}. The only difference is that we prefer to work with critical radii, while the computations in loc. cit. are mainly under the assumption of non-critical radii. Finally, we deduce that the bounded-by-$1$ functions $\mathcal{O}(\frX)^\circ$ on the character variety are all given by measures.

Let us write $\lVert \cdot \rVert_n$ for the Gau{\ss} norm on $B_n$ and consider for $s\in \RR$ the filtration
\begin{align*}
	F^s B_n:&=\{F\in B_n: \lVert F\rVert_n \leq p^{-s}\},\\
	F^{s,+} B_n:&=\{F\in B_n: \lVert F\rVert_n < p^{-s}\}.
\end{align*}
This filtration induces on $L$ the usual filtration by the absolute value. Let us write  $\gr^\bullet B_n$ and $\gr^\bullet L$ for the associated graded rings and $\sigma\colon B_n \to \gr^\bullet B_n$ for the principal symbol map, which is defined for $f\in F^sB_n\setminus F^{s,+} B_n$ as 
\[
\sigma(f):=f +F^{s,+} B_n\in F^{s} B_n/F^{s,+} B_n\subseteq \gr^\bullet B_n.
\]
 Let us write $\epsilon:=\sigma(p)\in \gr^\bullet L$. Then $\gr^\bullet L=k_L[\epsilon,\epsilon^{-1}]$ and
we have
\[
	\gr^\bullet B_n=\gr^\bullet L[X_1,X_2],
\]
with $X_i=\sigma(x_i)$, cf. \cite[Theorem 4.5]{SchneiderTeitelbaum2003}.

As explained above, the dominating terms of the logarithm series $\log(1+T)=\sum_{k\geq 1}(-1)^{k-1}\frac{T^k}{k} $ for the critical radius $\rho_n$ are $\frac{T^{p^n}}{p^n}$ and $\frac{T^{p^{n-1}}}{p^{n-1}}$. Hence
\[
	\sigma(H)=\pm \epsilon^{-n}\left[ \epsilon\left( X_2^{p^{n-1}}-\ol{v}_2X_1^{p^{n-1}} \right)+\left(X_2^{p^{n}}-\ol{v}_2X_1^{p^{n}}\right) \right].
\]
For later reference, and since it is crucial to the argument, we note that $\sigma(H)$ is irreducible:
\begin{lemma}
The element $\sigma(H)\in \gr^\bullet B_n$ is irreducible.
\end{lemma}
\begin{proof}
Since $\epsilon^{-n}\in \gr^\bullet B_n$ is a unit, it suffices to show that
\[
	\epsilon\left( X_2^{p^{n-1}}-\ol{v}_2X_1^{p^{n-1}} \right)+\left(X_2^{p^{n}}-\ol{v}_2X_1^{p^{n}}\right)\in k_L[\epsilon,\epsilon^{-1},X_1,X_2]
\]
is irreducible. We first prove that $X_2^{p^{n-1}}-\ol{v}_2X_1^{p^{n-1}}$ and $X_2^{p^{n}}-\ol{v}_2X_1^{p^{n}}$ are coprime in $k_L[X_1,X_2]$: We have
\[
	\left(X_2^{p^{n}}-\ol{v}_2X_1^{p^{n}}\right)-\left(X_2^{p^{n-1}}-\ol{v}_2X_1^{p^{n-1}}\right)^p=(\ol{v}_2^p-\ol{v}_2)X_1^{p^n}.
\]
Since $\ol{v}_2 \in k_L\setminus \mathbb{F}_p$ we have $(\ol{v}_2^p-\ol{v}_2)\neq 0$. We deduce that any common irreducible factor must be associated to $X_1$. However, neither $X_2^{p^{n-1}}-\ol{v}_2X_1^{p^{n-1}}$ and $X_2^{p^{n}}-\ol{v}_2X_1^{p^{n}}$ is divisible by $X_1$. This shows that $X_2^{p^{n-1}}-\ol{v}_2X_1^{p^{n-1}}$ and $X_2^{p^{n}}-\ol{v}_2X_1^{p^{n}}$ are co-prime in $k_L[X_1,X_2]$. Now view
\[
	\epsilon\left( X_2^{p^{n-1}}-\ol{v}_2X_1^{p^{n-1}} \right)+\left(X_2^{p^{n}}-\ol{v}_2X_1^{p^{n}}\right)\in k_L[X_1,X_2][\epsilon]
\]
as a degree $1$ polynomial in $\epsilon$ over $k_L[X_1,X_2]$. Its coefficients are co-prime, hence Gau{\ss}'s Lemma implies that it is irreducible in $k_L[\epsilon,X_1,X_2]$. Since $\epsilon\left( X_2^{p^{n-1}}-\ol{v}_2X_1^{p^{n-1}} \right)+\left(X_2^{p^{n}}-\ol{v}_2X_1^{p^{n}}\right)$ is not associated to $\epsilon$, localization at $\epsilon$ preserves irreducibility.
\end{proof}

By the above Lemma, $(\sigma(H))\subseteq \gr^\bullet B_n$ is a prime ideal. Let us now equip $A_n$ with the quotient norm
\[
	\lVert f\rVert_{n,\quot}:=\inf_{\substack{F\in B_n \\ F\bmod{(H)}=f}}\lVert F\rVert_n.
\]
Here, note that the ideal $(H)$ is closed, so the quotient seminorm is actually a norm.
Since the (additive) value group of the Gau{\ss} norm $\lVert\cdot \rVert_n$ on $B_n$ is discrete, we can find for any non-zero $f\in A_n$ a lift $F\in B_n$ with $\lVert F\rVert_n=\lVert f\rVert_{n,\quot}$.
\begin{lemma}\label{lem:char-minimal-lifts}
	Take $0\neq f\in A_n$ and let $F\in B_n$ be a lift of $f$. The following are equivalent:
	\begin{enumerate}
	\item $\lVert F\rVert_n=\lVert f\rVert_{n,\quot}$,
	\item $\sigma(F)\notin (\sigma(H))$.
	\end{enumerate}
\end{lemma}
\begin{proof}
For any lift, we have $\lVert F\rVert_n\geq \lVert f\rVert_{n,\quot}$. Suppose we have $\lVert F\rVert_n\neq \lVert f\rVert_{n,\quot}$, then $\lVert F\rVert_n> \lVert f\rVert_{n,\quot}$. Since we can always find a representative in $F+(H)$ of minimal norm $\lVert f\rVert_{n,\quot}$, there exists $G\in B_n$ such that $\lVert F-GH\rVert_n=\lVert f\rVert_{n,\quot}< \lVert F\rVert_n$. This implies that the dominating terms of $F$ and $GH$ cancel, i.e.,
\[
	\sigma(F)=\sigma(GH)=\sigma(G)\sigma(H) \in (\sigma(H)).
\]
Conversely, suppose we have $\sigma(F)\in (\sigma(H))$. Then we can find $G\in B_n$ such that $\sigma(F)=\sigma(H)\sigma(G)$. This implies $\lVert F-GH \rVert_n<\lVert F\rVert_n$. Hence we found a lift of strictly smaller norm and deduce $\lVert F\rVert_n\neq \lVert f\rVert_{n,\quot}$.
\end{proof}

\begin{corollary}\label{cor:quot-equals-spc}
We have $\lVert \cdot \rVert_{n,\quot}=\lVert \cdot \rVert_{n,\spc}$ on $A_n$.
\end{corollary}
\begin{proof}
According to \cite[6.2.3., Proposition 3]{BGR1984} it suffices to show that $\lVert \cdot \rVert_{n,\quot}$ is (power-) multiplicative. For non-zero $f,g\in A_n$, choose minimal lifts $F,G\in B_n$, i.e. $\lVert F\rVert_n=\lVert f\rVert_{n,\quot}$ and $\lVert G\rVert_n=\lVert g\rVert_{n,\quot}$. Lemma \ref{lem:char-minimal-lifts} implies $\sigma(F)\notin (\sigma(H))$ and $\sigma(G)\notin (\sigma(H))$. Since $(\sigma(H))$ is a prime ideal, we conclude $\sigma(FG)=\sigma(F)\sigma(G)\notin (\sigma(H))$. Again by Lemma \ref{lem:char-minimal-lifts}, we get $\lVert FG\rVert_n=\lVert fg\rVert_{n,\quot}$, i.e. we have
\[
	\lVert f\rVert_{n,\quot} \cdot \lVert g\rVert_{n,\quot}=\lVert F \rVert_n \cdot \lVert G \rVert_n = \lVert FG \rVert_n=\lVert fg\rVert_{n,\quot}. \qedhere
\]
\end{proof}

\begin{proof}[Proof of Theorem \ref{thm:main} in the unramified case]
Let us first recall (Lemma 1.15 of \cite{BSX}) that the map $o_L\llbracket o_L\rrbracket\to \mathcal{O}(\frX)$ is injective and has image in $\mathcal{O}(\frX)^\circ$. In the following, we prove that the image is all of $\mathcal{O}(\frX)^\circ$.
For $n$ a positive integer, write $A_n$ for the ring of functions on the critical radius $\rho_n$ as before. We have
\[
	\mathcal{O}(\frX)^\circ=\varprojlim_n A_n^\circ, 
\]
with $A_n^\circ=\{f\in A_n: \lVert f\rVert_{n,\spc}\leq 1\}$. Let us write $\Lambda_n$ for the image of $\Lambda:=o_L\llbracket o_L \rrbracket$ in $A_n^\circ$. Fix now $n\geq 1$ and consider $N>n$. We claim that $\im(A_N^\circ \to A_n^\circ) \subseteq \Lambda_n+ p^{p^{N-n}-1}A_n^\circ$.\\
\textit{Proof of the Claim:} Take $f\in A_N^\circ$. Choose a lift $F\in B_N$ with $\lVert F \rVert_N=\lVert f \rVert_{N,\quot}$. By Corollary \ref{cor:quot-equals-spc}, we have $\lVert f \rVert_{N,\quot}=\lVert f \rVert_{N,\spc}$. Hence $F\in B_N^\circ$. If we write $F=\sum_{\ul{k}\in \NN_0^2} c_{\ul{k}} \ul{x}^{\ul{k}}$ (with the notation $\ul{x}^{\ul{k}}:=x_1^{k_1}x_2^{k_2}$), this implies that for all $\ul{k}\in \NN_0^2$ we have $|c_{\ul{k}}|_p\leq \rho_N^{-|\ul{k}|}$.

We can approximate $F$ by the polynomial
\[
	P_N(x_1,x_2):=\sum_{|\ul{k}|<p^N-p^{N-1}} c_{\ul{k}} \ul{x}^{\ul{k}}\in L[x_1,x_2].
\]
Note that $P_N(x_1,x_2)$ actually has integral coefficients because for $|\ul{k}|<p^N-p^{N-1}$, we have
\[
	|c_{\ul{k}}|_p\leq \rho_N^{-|k|}< \rho_N^{-(p^N-p^{N-1})}=p.
\]
Since $L$ has (additive) value group $\ZZ$, this implies $c_{\ul{k}} \in o_L$. Hence, we have
\[
	P_N(x_1,x_2)\in o_L[x_1,x_2].
\]
We define
\[
	R_N:=F-P_N=\sum_{|\ul{k}|\geq p^N-p^{N-1}} c_{\ul{k}} \ul{x}^{\ul{k}} \in B^\circ_N.
\]
Using restriction, let us view $R_N$ as an element of $ B_n$:
\[
	\lVert R_N \rVert_n =\sup_{\ul{k}} |c_{\ul{k}}| \left(\frac{\rho_n}{\rho_N}\right)^{|\ul{k}|} \rho_N^{|\ul{k}|}\leq \left(\frac{\rho_n}{\rho_N}\right)^{p^N-p^{N-1}} \lVert R_N\rVert_N\leq \left(\frac{\rho_n}{\rho_N}\right)^{p^N-p^{N-1}}.
\]
With $\frac{\rho_n}{\rho_N}= p^{-\frac{1}{p^n-p^{n-1}}+\frac{1}{p^N-p^{N-1}}}$ we get
\[
	\lVert R_N \rVert_n\leq \left(\frac{\rho_n}{\rho_N}\right)^{p^N-p^{N-1}}\leq p^{-p^{N-n}+1},
\]
i.e. $R_N\in p^{p^{N-n}-1}B_n^\circ$, i.e.
\[
	F\in o_L[x_1,x_2]+p^{p^{N-n}-1}B_n^\circ.
\]
Passing to the quotient $B_n^\circ \twoheadrightarrow A_n^\circ$ shows that $f\in \Lambda_n+ p^{p^{N-n}-1}A_n^\circ$. This finishes the proof of the claim.\\

Since $A_N$ is an integral domain for every $N\geq 1$ (see the above discussion), the restriction maps $A_N\to A_n$ for $N\geq n$ are injective, cf. \cite[Proposition 4.2]{ABB2018}. It follows that the induced map
\[
	\mathcal{O}(\frX)=\varprojlim_{N}A_N\hookrightarrow A_n
\]
is injective  for every $n\geq 1$. Given $f\in \mathcal{O}(\frX)^\circ$, write $f_n:=f|_{\frX_n}\in \mathcal{O}(\frX_n)^\circ=A_n^\circ$. By the claim, we have
\[
	f_n\in \bigcap_{N\geq n} \Lambda_n + p^{p^{N-n}-1}A_n^\circ=\Lambda_n.
\]
Here, the last equality follows since $\Lambda_n$ is closed as the image of the compact set $o_L\llbracket o_L\rrbracket$ under the continuous map $o_L\llbracket o_L\rrbracket\to A_n$.
We deduce $f\in o_L\llbracket o_L\rrbracket$ as the restriction map $\mathcal{O}(\frX)\hookrightarrow A_n$ is injective and $o_L\llbracket o_L\rrbracket$ maps isomorphically to $\Lambda_n$ under this restriction map.
\end{proof}

\section{The ramified case}
The proof in the ramified case is quite similar to the unramified case with one small modification. Let us sketch the proof in the ramified case. 
Let $L/\Qp$ be a ramified quadratic extension. If $\pi$ is a uniformizer of $L$, then
$v_1=1$ and $v_2=\pi$ generate $o_L$ as a $\Zp$-module. Before we proceed, let us briefly discuss why the above argument needs a modification. The character variety is again locally cut out by the equation
\[
	H(x_1,x_2):=\log(1+x_2)-\pi\log(1+x_1),
\]
however, for the previous choice of $B_n$ with radius $\rho_n=p^{-1/(p^n-p^{n-1})}$, we run into the following problem. Since $\sigma(\pi)$ is not of degree zero, we have 
\[
 \sigma(H)=\pm \epsilon^{-n} \left( \epsilon X_2^{p^{n-1}} +X_2^{p^{n}}\right),
\]
which is not prime, and for $n>1$ the quotient by the ideal $(\sigma(H))$ is not even reduced. However, the above proof can be fixed by considering two different radii. For $n\geq 1$, define
\[
	\alpha_n:=p^{-\frac{1}{p^n-p^{n-1}}-\frac{1}{4p^{n-1}}}, \quad \beta_n:=p^{-\frac{1}{p^n-p^{n-1}}-\frac{3}{4p^{n-1}}}.
\]
These families of radii satisfy
\begin{equation}\label{eq:radii-alpha-beta}
	p^{-1/2}\alpha_n^{p^{n-1}}=\beta_n^{p^{n-1}}.
\end{equation}
The importance of this equation (and hence why we have chosen $\alpha_n$ and $\beta_n$ as above) will become apparent in a second, when we compute the initial form of $H$.
Let us write
\[
	B_n:=L\langle \alpha_n^{-1}x_1 , \beta_n^{-1}x_2 \rangle
\]
and $A_n:=B_n/(H)$. The dominant term of $ \log(1+x_2)$ on $B_n$ is $\frac{x_2^{p^{n-1}}}{p^{n-1}}$, while the dominant term of $\pi \log(1+x_1)$ is $\pi \frac{x_1^{p^{n-1}}}{p^{n-1}}$. Note that both terms have the same Gau{\ss} norm because of \eqref{eq:radii-alpha-beta}. Let us write $\epsilon:=\sigma(\pi)$. Then $\gr^\bullet L=\mathbb{F}_p[\epsilon,\epsilon^{-1}]$ and
\[
	\gr^\bullet B_n =(\gr^\bullet L)[X_1,X_2]
\]
with $X_1=\sigma(x_1)$ and $X_2=\sigma(x_2)$. The ideal generated by the initial form of $H$ is
\[
	(\sigma(H))=\left( X_2^{p^{n-1}}-\epsilon X_1^{p^{n-1}} \right)\subseteq \gr^\bullet B_n.
\]
As in the unramified case, $X_2^{p^{n-1}}-\epsilon X_1^{p^{n-1}}$ is a polynomial of degree $1$ in $\epsilon$ with co-prime coefficients and hence irreducible.
This shows that $(\sigma(H))$ is a prime ideal which was the essential ingredient for proving that the spectral norm on $A_n$ coincides with the quotient norm. We can now use essentially the same argument as in the unramified case to prove the ramified case.

\begin{proof}[Proof of Theorem \ref{thm:main} in the ramified case]
Let $n$ be a positive integer and write $A_n$ for the ring of functions of radius $(\alpha_n,\beta_n)$ as above. We have
\[
	\mathcal{O}(\frX)^\circ=\varprojlim_n A_n^\circ, 
\]
with $A_n^\circ=\{f\in A_n: \lVert f\rVert_{n,\spc}\leq 1\}$. Let us write $\Lambda_n$ for the image of $\Lambda:=o_L\llbracket o_L \rrbracket$ in $A_n^\circ$. Fix now $n\geq 1$ and consider $N>n$.\\
\textit{Claim:} $\im(A_N^\circ \to A_n^\circ) \subseteq \Lambda_n+ p^{m_{N-n}}A_n^\circ$ with $m_{N-n}:=\left\lfloor \frac{(p+3)}{8} \cdot\left(p^{N-n-1}-\frac{1}{p}\right) \right\rfloor$.\\
\textit{Proof of the Claim:} Let $f\in A_N^\circ$. Choose a lift $F\in B_N$ with $\lVert F \rVert_N=\lVert f \rVert_{N,\quot}$. By the analogue of Corollary \ref{cor:quot-equals-spc}, we have $\lVert f \rVert_{N,\quot}=\lVert f \rVert_{N,\spc}$. Hence $F\in B_N^\circ$. If we write $F=\sum_{\ul{k}\in \NN_0^2} c_{\ul{k}} \ul{x}^{\ul{k}}$ (with the notation $\ul{x}^{\ul{k}}:=x_1^{k_1}x_2^{k_2}$), this implies for all $\ul{k}\in \NN_0^2$
\[
	|c_{\ul{k}}|_p\leq \alpha_N^{-k_1}\beta_N^{-k_2}.
\]
We can approximate $F$ by the polynomial
\[
	P_N(x_1,x_2):=\sum_{|k|< \frac{p^{N-1}-p^{N-2}}{2}} c_{\ul{k}} \ul{x}^{\ul{k}}\in L[x_1,x_2].
\]
Note that $P_N(x_1,x_2)$ has actually integral coefficients because for $|\ul{k}|<\frac{p^{N-1}-p^{N-2}}{2}$, we have
\[
	|c_{\ul{k}}|_p\leq \alpha_N^{-k_1}\beta_N^{-k_2}= p^{ \frac{k_1(4+p-1)}{4p^{N-1}(p-1)}+\frac{k_2(4+3(p-1)}{4p^{N-1}(p-1)}}\leq p^{ \frac{k_1+k_2}{p^{N-2}(p-1)}}<p^{1/2}.
\]
Since $L$ has (additive) value group $\frac{1}{2}\ZZ$, this implies $c_{\ul{k}}\in o_L$. Hence, we have
\[
	P_N(x_1,x_2)\in o_L[x_1,x_2].
\]
We define
\[
	R_N:=F-P_N=\sum_{|\ul{k}|\geq \frac{p^{N-1}-p^{N-2}}{2}} c_{\ul{k}} \ul{x}^{\ul{k}} \in B^\circ_N.
\]
Using restriction, let us view $R_N$ as an element of $ B_n$:
\begin{align*}
	\lVert R_N \rVert_n &=\sup_{\ul{k}} |c_{\ul{k}}| \left(\frac{\alpha_n}{\alpha_N}\right)^{k_1} \left(\frac{\beta_n}{\beta_N}\right)^{k_2} \alpha_N^{k_1}\beta_N^{k_2}\\
	&\leq \left(\frac{\alpha_n}{\alpha_N}\right)^{\frac{p^{N-1}-p^{N-2}}{2}} \lVert R_N\rVert_N\\
	&\leq \left(\frac{\alpha_n}{\alpha_N}\right)^{\frac{p^{N-1}-p^{N-2}}{2}}
\end{align*}
Here, we have used $0<\frac{\beta_n}{\beta_N}\leq \frac{\alpha_n}{\alpha_N}<1$ for $n<N$ in the first inequality.
Let us write $\log_p$ for the logarithm with base $p$. We have
\begin{align*}
\log_p\left(\frac{\alpha_n}{\alpha_N}\right)&= {-\frac{1}{p^n-p^{n-1}}-\frac{1}{4p^{n-1}}+\frac{1}{p^N-p^{N-1}}+\frac{1}{4p^{N-1}}}\\
&= -\frac{(p+3)}{4} \cdot\frac{p^{N-n-1}-\frac{1}{p}}{p^{N-1}-p^{N-2}}
\end{align*}
we get
\begin{align*}
	\log_p\lVert R_N \rVert_n &\leq \frac{p^{N-1}-p^{N-2}}{2}\cdot \log_p \left(\frac{\alpha_n}{\alpha_N}\right) \\
	&=-\frac{(p+3)}{8} \cdot\left(p^{N-n-1}-\frac{1}{p}\right).
\end{align*}
i.e. we get
\[
	F\in o_L[x_1,x_2]+p^{m_{N-n}}B_n^\circ.
\]
Passing to the quotient $B_n^\circ \twoheadrightarrow A_n^\circ$ shows that $f\in \Lambda_n+ p^{m_{N-n}}A_n^\circ$. This finishes the proof of the claim.\\

Since $m_{N-n}\to \infty$ as $N\to \infty$, one can argue as in the unramified case to deduce that $o_L\llbracket o_L\rrbracket\to \mathcal{O}(\frX)^\circ$ is an isomorphism.
\end{proof}

\section{Associated orders}
\label{assord}

We now prove Theorem \ref{ordcrit:intro}. Let $L = \QQ_{p^2}$ and $\pi=p$, and let $L_n$ be the extension of $L$ generated by the $p^n$-torsion points of the Lubin--Tate group attached to $L$ and $p$. Let $G_n = \Gal(L_n/L)$ and let $\Gamma_L = \varprojlim G_n$. Let $\calA_n = \{ f \in L[G_n]$ such that $f(o_n) \subset o_n \}$. There are restriction maps $\pi_n : \calA_{n+m} \to \calA_n$ for all $m,n \geq 0$.

\begin{theorem}
\label{ordcrit}
For all $n \geq 0$, there exists $m$ such that $\pi_n(\calA_{n+m}) = o_L[G_n]$.
\end{theorem}

\begin{proof}
Let $o_\infty $ be the $p$-adic completion of $\cup_{n \geq 1} o_n$. By Corollary 5.1.7 of \cite{ArdakovBerger2024}, we have $\Lambda(X) = o_L \dcroc{o_L}$ if and only if every $\Gamma_L$-equivariant $o_L$-linear map $o_\infty \to o_\infty$ comes from an element of $o_L \dcroc{\Gamma_L}$. We now prove that this condition implies the theorem. 

Since $L_n/L$ is abelian, the normal basis theorem implies that $L[G_n]$ is the set of $\Gamma_L$-equivariant $L$-linear maps $L_n \to L_n$.  If $f \in \End_{o_L}^{\Gamma_L}(o_\infty)$ then $f(o_n) \subset o_n$ by the Ax-Sen-Tate theorem. The image of such an $f$ in $L[G_n]$ therefore lies in $\calA_n$. This gives rise to an isomorphism $\End_{o_L}^{\Gamma_L}(o_\infty) = \varprojlim \calA_n$. Suppose that $\End_{o_L}^{\Gamma_L}(o_\infty) = o_L \dcroc{\Gamma_L}$. Let $\calM_n = o_L [G_n] \subset \calA_n$ and let $\calQ_n = \calA_n / \calM_n$ which is a finite group. We have an exact sequence $0 \to \calM_n \to \calA_n \to \calQ_n \to 0$. The maps $\pi_n : \calM_{n+m} \to \calM_n$ are surjective, so that the Mittag--Leffler condition holds and the sequence $0 \to \varprojlim \calM_n \to \varprojlim \calA_n \to \varprojlim \calQ_n \to 0$ is exact. If the map $\varprojlim \calM_n \to \varprojlim \calA_n$ is an isomorphism, then $\varprojlim \calQ_n = \{0\}$. Since each $\calQ_n$ is a finite group, this implies that for all $n \geq 0$, there exists $m$ such that the image of $\calQ_{m+n}$ in $\calQ_n$ is zero.
\end{proof}

\section{Concluding remarks}\label{sec:concluding-remarks}
For $d:=[L:\Qp]>2$, a direct generalization of the above arguments is not possible. Let us explain this in the case where $L$ is an unramified extension of degree $3$. Let us choose a $\mathbb{Z}_p$-basis $v_1=1,v_2,v_3$ of $o_L$. We keep the notation introduced in section \ref{sec:notation} with the choice of radii 
\[
	\rho_n:=p^{-\frac{1}{p^n-p^{n-1}}}.
\]
In this case, we have
\[
	A_n=B_n/I
\]
with $I=(H_2,H_3)$, where
\[
	H_2=\log(1+x_2)-v_2 \log(1+x_1), \quad H_3=\log(1+x_3)-v_3 \log(1+x_1).
\]
The principal symbols $\sigma(H_i)\in \gr^\bullet B_n=k_L[\epsilon,\epsilon^{-1},X_1,X_2,X_3]$ are given by
\begin{align*}
	\sigma(H_2)&=\pm \epsilon^{-n}\left[ \epsilon\left( X_2^{p^{n-1}}-\ol{v}_2X_1^{p^{n-1}} \right)+\left(X_2^{p^{n}}-\ol{v}_2X_1^{p^{n}}\right) \right]\\
	\sigma(H_3)&=\pm \epsilon^{-n}\left[ \epsilon\left( X_3^{p^{n-1}}-\ol{v}_3X_1^{p^{n-1}} \right)+\left(X_3^{p^{n}}-\ol{v}_3X_1^{p^{n}}\right) \right].
\end{align*}
Let us write $\sigma(I)$ for the ideal generated by all principal symbols of all elements of $I$. Since $\sigma(H_2),\sigma(H_3)$ forms a regular sequence in $\gr^\bullet B_n$, we have by \cite[Ch. III, Proposition 3.10]{NvO1979} that
\[
	\sigma(I)=(\sigma(H_2),\sigma(H_3)).
\]
Choose elements $\overline{w}_2,\overline{w}_3\in k_L$ such that $\overline{w}_2^{p^{n-1}}=\overline{v}_2$ and $\overline{w}_3^{p^{n-1}}=\overline{v}_3$, and define for $i=2,3$:
\begin{align*}
	L_i:&=X_i-\overline{w}_iX_1\\
	K_i:&=X_i^p-\overline{w}_iX_1^p.
\end{align*}
Then
\[
	\sigma(I)=\left( \epsilon L_2^{p^{n-1}}+K_2^{p^{n-1}}, \epsilon L_3^{p^{n-1}}+K_3^{p^{n-1}} \right).
\]
For $n\geq 2$, the element
\[
	L_3 K_2-K_3L_2 \in \gr^\bullet B_n
\]
is not contained in $\sigma(I)$, but
\[
	\left(L_3 K_2-K_3L_2\right)^{p^{n-1}}=L_3^{p^{n-1}}K_2^{p^{n-1}}-K_3^{p^{n-1}}L_2^{p^{n-1}}\in \sigma(I).
\]
Note that this implies that the quotient norm on $A_n$ is not multiplicative. Indeed, choose any $F\in B_n$ with
\[
	\sigma(F)=L_3 K_2-K_3L_2 \notin \sigma(I).
\]
Write $f\in A_n$ for the image of $F$ under the quotient map.
The same argument as in Lemma \ref{lem:char-minimal-lifts} shows that $\lVert F\rVert_n=\lVert f\rVert_{n,q}$. On the other hand, we have
\[
	\sigma(F^{p^{n-1}})\in \sigma(I),
\]
which implies $\lVert f^{p^{n-1}}\rVert_{n,q}< \lVert f\rVert_{n,q}^{p^{n-1}}$. This shows that the quotient norm is not power-multiplicative and hence it does not coincide with the spectral norm.

\section{AI usage}
The article was entirely written by the authors. All of the key ideas entering the proof were of a human nature (except for \S \ref{assord}). ChatGPT was used for proofreading and checking mathematical correctness. The specific choices of radii ($\alpha_n$ and $\beta_n$) in the ramified case were suggested by ChatGPT in the process of proofreading of an early version of the manuscript.

\bibliographystyle{amsalpha}
\bibliography{BFReferences}
\end{document}